\documentstyle[12pt]{article}

\newtheorem{proposition}{Proposition}[section]
\newtheorem{theorem}[proposition]{Theorem}

\newcommand{\PP}{\mbox{\sf I\hspace*{-.12em}P}}
\newcommand{\PG}[2]{\mbox{$\mbox{{\rm PG}}(#1,#2)$}}
\newcommand{\call}[1]{{\cal #1}}

\newcommand{\rist}[1]{\mbox{\raisebox{-1ex}{$|#1$}}}
\newcommand{\ml}{linear morphism}
\newcommand{\mli}{linear morphisms}
\newcommand{\applicaz}[3]{\mbox{$#1\,:\,#2
  \rightarrow#3$}}
\newcommand{\ach}{\alpha\chi}
\newcommand{\gfi}{\gamma\varphi}
\newcommand{\rad}[1]{\mbox{{\rm rad}}_{#1}}

\title{On linear morphisms of product spaces\thanks{This work was performed
in the context of the fourth protocol of
  scientific and technological cooperation between Italy and Austria
  (plan n.~10).
  Partial support was provided by the project ``Strutture Geometriche,
  Combinatoria e loro Applicazioni'' of M.U.R.S.T.}}

\author{Alessandro Bichara, Hans Havlicek, Corrado Zanella}

\date{October 28, 1999}

\begin{document}

\maketitle

\begin{abstract}
  Let $\chi$ be a \ml\ of the product of two projective spaces
  \PG nF\ and \PG mF\ into a projective space.
  Let $\gamma$ be the Segre embedding of such a product.
  In this paper we give some sufficient conditions for the existence
  of an automorphism $\alpha$ of the product space and a
  \ml\ of projective spaces $\varphi$, such that $\gfi=\ach$.\\
  A.M.S. classification number: 51M35. Keywords: Segre variety --
  product space -- projective embedding.
\end{abstract}


\section{Introduction}

A {\em semilinear space\/} is a pair $\Sigma=(\call P,\call G)$, where
$\call P$ is a set, whose elements are called {\em points\/}, and
$\call G\subset2^{\call P}$.
(In this paper ``$A\subset B$'' just means that $x\in A$ implies $x\in B$.)
The elements of $\call G$ are {\em lines\/}.
The axioms defining a semilinear space are the following:
({\em i\/})~$|g|\geq2$ for every line $g$;
({\em ii\/})~$\cup_{g\in\call G}g=\call P$;
({\em iii\/})~$g,h\in\call G$, $g\neq h$ $\Rightarrow$ $|g\cap h|\leq1$.
Two points $X,Y\in\call P$ are {\em collinear\/}, $X\sim Y$, if a line $g$
exists such that $X,Y\in g$ (for $X\neq Y$ we will also write $XY:=g$).
An {\em isomorphism\/} between the semilinear spaces $(\call P,\call G)$
and $(\call P',\call G')$ is a bijection \applicaz{\alpha}{\call P}{\call P'}\
such that both $\alpha$ and $\alpha^{-1}$ map lines onto lines.

The {\em join\/} of $\call M_1,\call M_2\subset \call P$ is:
\[
  \call M_1\vee\call M_2:=\call M_1\cup\call M_2\cup
  \bigg(
    \bigcup_{
        X_i\in\call M_i\atop X_1\sim X_2,\,X_1\neq X_2
    }X_1X_2
  \bigg).
\]

If $X$ is a point, we will often write $X$ instead of $\{X\}$.
Let $\PP'=(\call P',\call G')$ be a projective space.
A {\em \ml\/} \applicaz{\chi}{\Sigma}{\PP'}\ is a mapping of a subset
${\bf D}(\chi)$ of $\call P$ into $\call P'$ satisfying the following
axioms (L1) and (L2)~\cite{Bra73,Hav81}.
Here $X\chi=\emptyset$ for $X\in{\bf A}(\chi):=\call P\setminus{\bf D}(\chi)$.

(L1)\hspace{1cm}$(X\vee Y)\chi=X\chi\vee Y\chi$ for $X,Y\in\call P$,
$X\sim Y$;

(L2)\hspace{1cm}$X,Y\in\call P$, $X\chi=Y\chi$,  $X\neq Y$, $X\sim Y$
$\Rightarrow$ $\exists A\in XY$ such that $A\chi=\emptyset$.

${\bf D}(\chi)$ is the {\em domain} of $\chi$ and
${\bf A}(\chi)$ is the {\em exceptional set\/}.
The \ml\ $\chi$ is said to be {\em global\/} when
${\bf D}(\chi)=\call P$; is called {\em embedding\/} if it is global and
injective.
It should be noted the last definition is somewhat particular, since
for instance the inclusion of an affine space into its projective extension
is not an embedding.

Throughout this paper we will deal with the following:
a commutative field $F$; two natural numbers $n\geq 2$ and $m\geq1$;
$\PP_1=\PG nF=(\call P_1,\call G_1)$, that is
the projective space of dimension $n$ over $F$;
$\PP_2=\PG mF=(\call P_2,\call G_2)$;
$\overline{\PP}=\PG{nm+n+m}F=(\overline{\call P},\overline{\call G})$;
a further projective space $\PP'=(\call P',\call G')$.

The {\em product\/} of $\PP_1$ and $\PP_2$ is the
semilinear space $\PP_1\times\PP_2=(\call P^*,\call G^*)$, where
$\call P^*=\call P_1\times\call P_2$, and the elements of $\call G^*$ are of
two kinds: $ X_1 \times g_2$ with $X_1\in\call P_1$, $g_2\in\call G_2$, and
$g_1\times X_2 $ with $g_1\in\call G_1$, $X_2\in\call P_2$.

We shall be concerned with a \ml\
\applicaz{\chi}{\PP_1\times\PP_2}{\PP'}, and a {\em regular\/}
embedding \applicaz{\gamma}{\PP_1\times\PP_2}{\overline{\PP}}.
The word ``regular'' means that the projective closure
$[(\call P_1\times\call P_2)\gamma]$ of $(\call P_1\times\call P_2)\gamma$
has dimension $mn+m+n$
(which is the greatest possible one, see~\cite{Zan96}).
Our purpose is to generalize the main result of~\cite{Zan96}, which is the
following: if $\chi$ is an embedding, then there exist an automorphism
$\alpha'$ of
$\PP_2$ and a \ml\ \applicaz{\varphi}{\overline{\PP}}{\PP'}\ such that
$\gfi=(\mbox{id}_{\call P_1},\alpha')\chi$.

We define the {\em first\/} and {\em second radical\/} of $\chi$ by
\[\begin{array}{l}
  \rad1\chi:=\{X\in\call P_1|X\times\call P_2\subset{\bf A}(\chi)\};\\
  \rad2\chi:=\{Y\in\call P_2|\call P_1\times Y\subset{\bf A}(\chi)\}.
\end{array}\]
$\chi$ is a {\em degenerate\/} \ml\ if
$\rad1\chi\cup\rad2\chi\neq\emptyset$.

Clearly each radical $\rad i\chi$ is a subspace of $\PP_i$, $i=1,2$.
The following proposition can be useful in reducing the investigation
to nondegenerate \mli.

\begin{proposition}\label{nondegen}
  Let $\rad i\chi$ and $D_i$ be complementary subspaces of
  $\PP_i$ ($i=1,2$).
  Let \applicaz{\pi_i}{\PP_i}{D_i}\ be the projection onto $D_i$ from
  $\rad i\chi$, and $\chi':=\chi\rist{D_1\times D_2}$.
  Then, for every $(X,Y)\in\call P_1\times\call P_2$, it holds:
  \begin{equation}
    \{(X,Y)\}\chi=(X\pi_1\times Y\pi_2)\chi'.
  \end{equation}
\end{proposition}

{\em Proof\/}.
If $X\in\rad 1\chi$ or $Y\in\rad2\chi$ the statement is trivial.
Otherwise, let $X'\in\rad1\chi$ such that $X$, $X'$ and $X\pi_1$ are collinear.
We have
$\{(X,Y)\}\chi\in\{(X',Y)\}\chi\vee\{(X\pi_1,Y)\}\chi=\{(X\pi_1,Y)\}\chi$;
since $X\neq X'$, the set $\{(X,Y)\}\chi$ is empty if and only if
$\{(X\pi_1,Y)\}\chi$ is.
Then
$\{(X,Y)\}\chi=(X\pi_1\times Y)\chi$.
A similar argument on $Y$ yields the result.$\Box$

We prove next the following result which will be useful later on.
\begin{proposition}\label{slm}
  Let $\ell_1$, $\ell_2$ be two lines, \applicaz{\mu}{\ell_1\times\ell_2}{\PP'}
  a \ml, and
  $A$ a point of $\ell_1$ such that $A\times \ell_2$
  is contained in ${\bf D}(\mu)$.
  Then the exceptional set  ${\bf A}(\mu)$ is either a line, or
  $|{\bf A}(\mu)|\leq2$.
  Furthermore:\\
  (i)~If $|{\bf A}(\mu)|=1$ and $(X_1,P)$ is the only exceptional point,
  then, for $X\in\ell_1\setminus\{A,X_1\}$, one of the following possibilities
  occurs:\\
  (a)~$(A\times \ell_2)\mu=(X\times \ell_2)\mu$; in this case
  $\{(A,P)\mu\}=(X_1\times \ell_2)\mu$;\\
  (b)~$(A\times \ell_2)\mu\cap(X\times \ell_2)\mu$ is the point
  $(A,P)\mu$.\\
  (ii) If $|{\bf A}(\mu)|=2$, and $(X_1,P)$, $(X'_1,P')$ are
  the exceptional points, then $X_1\neq X'_1$, $P\neq P'$,
  $(X'_1\times \ell_2)\mu=\{(A,P)\mu\}$, and
  $(\ell_1\times \ell_2)\mu=(A\times \ell_2)\mu$.
\end{proposition}

{\em Proof\/}.
If ${\bf A}(\mu)$ contains a line, then such a line is of type
$P\times \ell_2$, $P\in\ell_1$.
By the assumption $A\times \ell_2\subset{\bf D}(\mu)$, there exist no
more exceptional points.

Now assume that ${\bf A}(\mu)$ does not contain lines, and that
$(X_1,P),(X'_1,P')\in{\bf A}(\mu)$.
This implies $X_1\neq X'_1$, $P\neq P'$.
Next,
\begin{equation}\label{eslm1}
  (\ell_1\times P)\mu=\{(A,P)\mu\}\neq\{(A,P')\mu\}=(\ell_1\times P')\mu,
\end{equation}
\begin{equation}\label{eslm2}
  \begin{array}{l}
    \{(A,P)\mu\}=\{(X'_1,P)\mu\}=(X'_1\times \ell_2)\mu,\\
    \{(A,P')\mu\}=\{(X_1,P')\mu\}=(X_1\times \ell_2)\mu.
  \end{array}
\end{equation}

Let $g$ be a line of $\ell_1\times\ell_2$.
If $g$ is of type $Y \times \ell_2$, then $g$ intersects $\ell_1\times P$
and $\ell_1\times P'$, hence $g\mu\subset(A\times \ell_2)\mu$ follows from
(\ref{eslm1}).
Otherwise, the same conclusion follows from $(\ref{eslm2})$.

Now assume that a further point $(X''_1,P'')\in{\bf A}(\mu)$ exists.
We have
$P\neq P''\neq P'$, and
\[
  (X''_1 \times \ell_2)\mu=\{(X''_1,P)\mu\}=\{(X''_1,P')\mu\}.
\]
This contradicts $(X''_1,P)\mu=(A,P)\mu\neq(A,P')\mu=(X''_1,P')\mu$.

Now we have to prove ({\em i\/}).
By assumption, $(X_1 \times \ell_2)\mu$ is a point $Z$.

If case ({\em a\/}) occurs, take a $Q\in\ell_2\setminus P$.
Since $(\ell_1\times Q)\subset{\bf D}(\mu)$ and $(X_1,Q)\mu=Z$,
it holds $(X,Q)\mu\neq Z$.
Furthermore,
$(\ell_1\times Q )\mu=( A \times \ell_2)\mu$ implies
$Z\in(A\times \ell_2 )\mu$ and $(A,P)\mu=Z$.

If case ({\em a\/}) does not occur, then
$(A,P)\mu=(X,P)\mu$ implies ({\em b\/}).$\Box$


\section{The first decomposition of $\chi$}

In this section we will prove a universal property of
a regular embedding $\gamma$
with respect to a given class of \mli.
We assume that the \ml\ $\chi$ satisfies the following properties:\\
({\em i\/}) There exist a basis $\call B=\{A_0,A_1,\ldots,A_m\}$ of $\PP_2$
and a plane $\call E\subset\call P_1$, such that for every
$i=1,\ldots,m$, $\chi\rist{\call E\times\{A_0,A_i\}}$ is global and has rank%
\footnote{The rank is the projective dimension of the projective closure
of $(\call E\times\{A_0,A_i\})\chi$.}
greater than 2.\\
({\em ii\/}) There exists $A\in\call P_1$ such that
$\dim( A \times\call P_2)\chi=m$.

\begin{proposition}\label{proiett}
  Assume $g\in\call G_1$, $P^*_1,P^*_2\in\call P_2$, $P^*_1\neq P^*_2$, and
  \begin{equation}\label{eproiett}
    g_i:=(g\times P^*_i )\gamma,\hspace{2em} i=1,2.
  \end{equation}
  Then the mapping
  \applicaz f{g_1}{g_2}, determined by
  $(P,P^*_1)\gamma f:=(P,P^*_2)\gamma$ for every $P\in g$,
  is a projectivity.
\end{proposition}

{\em Proof\/}.
Let $P^*_0\in P^*_1P^*_2\setminus\{P^*_1,P^*_2\}$ and
$g_0:=(g\times P^*_0 )\gamma$.
For every $X\in g_1$ it holds
$\{Xf\}=( X \vee g_0)\cap g_2$.$\Box$

\begin{proposition}\label{generano}
  Let $g\in\call G_1$, $g'\in\call G_2$, $P^*_0,P^*_1,P^*_2$
  three distinct points on $g'$, and
  \begin{equation}\label{egenerano}
    h_i:=(g\times P^*_i )\chi,\hspace{2em}  i=0,1,2.
  \end{equation}
  Then, for every permutation $(i,j,k)$ of $(0,1,2)$, it holds
  $h_i\subset h_j\vee h_k$.$\Box$
\end{proposition}

\begin{proposition}\label{perv4}
  With the assumptions of props.~\ref{proiett} and \ref{generano},
  if $h_1$ and $h_2$ are lines of $\PP'$ and
  \applicaz{\sigma'}{h_1}{h_2} is the mapping
  \begin{equation}\label{eperv4}
    (P,P_1^*)\chi\sigma':=(P,P^*_2)\chi,\hspace{2em}P\in g,
  \end{equation}
  then \applicaz{\sigma}{g_1}{g_2}, defined by
  \begin{equation}\label{defsigma}
    \sigma:=((\gamma^{-1}\chi)\rist{g_1})\sigma'((\gamma^{-1}\chi)\rist{g_2})^{-1},
  \end{equation}
  is a projectivity.
\end{proposition}

{\em Proof\/}.
Indeed $\sigma$ coincides with the mapping $f$ of prop.~\ref{proiett}.$\Box$

\begin{proposition}\label{esistefi}
  $(\gamma^{-1}\chi)\rist{(\call P_1\times\call B)\gamma}$
  can be extended to a \ml\\
  \applicaz\varphi{\overline{\PP}}{\PP'}.
\end{proposition}

{\em Proof\/}.
We show the following statement by induction on $t=0,\ldots,m$:\\
{\em $(S)_t$\hspace{5mm}There exists a \ml\
\applicaz{\varphi_t}{\left[(\call P_1\times\{A_0,A_1,\ldots,A_t\})\gamma
\right]}{\PP'}\\
which extends
$(\gamma^{-1}\chi)\rist{(\call P_1\times\{A_0,A_1,\ldots,A_t\})\gamma}$.}\\
For $t=0$ there is nothing to extend.
Suppose $(S)_t$ is true for a $t<m$.
Let $P^*_1=A_0$,
$P^*_2=A_{t+1}$, take
$P^*_0\in P^*_1P^*_2\setminus\{P^*_1,P^*_2\}$
so, that the dimension of
$U^*:=(\call E\times\{P^*_0\})\chi$ is minimal.
Define:
\[
  \psi:=\gamma^{-1}\chi\rist{(\call P_1\times P^*_2 )\gamma}.
\]
We use theor.~1.6 in \cite{Hav81}, case (V4),
in order to extend the pair of \mli\
$\varphi_t$ and $\psi$ to a \ml
\[
\applicaz{\varphi_{t+1}}{\left[(\call P_1\times\{A_0,A_1,\ldots,A_{t+1}\})
\gamma\right]}{\PP'}.
\]
By the regularity of $\gamma$, it holds
$((\call P_1\times P^*_2 )\gamma)\cap\left[(\call P_1\times\{A_0,A_1,\ldots,
A_t\})\gamma\right]=\emptyset$.
Next, it is enough to prove the existence of a line
$g\subset\call E$ such that the mapping
$\sigma'$ defined in (\ref{eperv4}) is a projectivity: indeed, in this case,
the mapping in  (\ref{defsigma}) is
$(\varphi_t\rist{g_1})\sigma'(\psi\rist{g_2})^{-1}$.

{\em Case 1: There exists a line $g\subset\call E$ such that the lines
$h_1$ and $h_2$ (cf.~(\ref{egenerano})) are skew.}\\
By prop.~\ref{generano}, any two of $h_0$, $h_1$ and $h_2$ are skew lines.
In this case, the same argument as in prop.~\ref{proiett} shows that
the mapping \applicaz{\sigma'}{h_1}{h_2},
which is defined in (\ref{eperv4}), is a projectivity.

{\em Case 2: $\dim U^*\leq0$.}\\
There exists a line, say $g$, lying on $\call E$ and such that
$g\times P^*_0 \subset{\bf A}(\chi)$.
For every $P\in g$,
\[
  (P,P^*_1)\chi\in(P,P^*_0)\chi\vee(P,P^*_2)\chi=\{(P,P^*_2)\chi\}.
\]
So, the mapping in (\ref{eperv4}) is the identity.

{\em Case 3: $\dim U^*=1$.}\\
In this case ${\bf A}(\chi)\cap(\call E\times P^*_0 )$
is a point $(Q,P^*_0)$ of $\PP_1\times\PP_2$.
If every line $h$ through $Q$ and in $\call E$ satisfies
$(h\times P^*_1 )\chi=(h\times P^*_2 )\chi$, then
$(\call E\times P^*_1 )\chi=(\call E\times P^*_2 )\chi$, and this contradicts
the assumptions.
Thus there is a line of $\call E$ and through $Q$, say $g$, such that,
according to definitions
(\ref{egenerano}), $h_1\neq h_2$.
The sets $h_1$, $h_2$ are lines and  $h_0$ is a point not belonging to
$h_1\cup h_2$ (cf.\ prop.~\ref{generano}).
For every $P\in g$ it holds
$(P,P^*_2)\chi\in h_0\vee\{(P,P^*_1)\chi\}$; therefore the mapping in
(\ref{eperv4}) is the perspectivity between $h_1$ and $h_2$ with center
$h_0$.

{\em Case 4: $\dim U^*=2$.}\\
Assume that we are not in case 1.
Then for every line $h$ of $\call E$ we have
$(h\times P^*_1 )\chi\cap(h\times P^*_2 )\chi\neq\emptyset$.
There is a line $g$ in $\call E$ such that
$\left((g\times P^*_1 )\chi\right)\cap
\left((\call E\times P^*_2 )\chi\right)$
is a point $P^*$.
It holds:
\begin{equation}\label{ecaso4}
  P^*=(P_1,P^*_1)\chi=(P_2,P^*_2)\chi\hspace{2em}{\rm with}\ P_1,P_2\in g.
\end{equation}
There exists a line $g'$ in $\call E$ such that $P_1\in g'$, $g'\neq g$ and
$(g'\times P^*_1 )\chi\cap((\call E\times P^*_2 )\chi)=P^*$.
Next, there is a point $P_2'\in g'$ such that
$(P_1,P^*_1)\chi=(P_2',P^*_2)\chi$.
It follows $P'_2=P_2=P_1$.
By (\ref{ecaso4}), there is $P'_0\in P^*_1P^*_2$ such that
$(P_1,P'_0)\in{\bf A}(\chi)$.
But then the dimension of $U^*$ is not minimal, a contradiction.$\Box$

Now we use the assumption ({\em ii\/}) at the beginning of this section:
\begin{equation}\label{evert}
  \begin{array}{ccl}
    ( A \times\call P_2)\chi&=&(A,A_0)\chi\vee(A,A_1)\chi\vee\ldots\vee
    (A,A_m)\chi=\\
    &=&(A,A_0)\gfi\vee(A,A_1)\gfi\vee\ldots\vee
    (A,A_m)\gfi=\\
    &=&( A \times\call P_2)\gfi.
  \end{array}
\end{equation}
Then for every $X\in\call P_2$ there is exactly one element of $\call P_2$, say
$Y$, such that
\begin{equation}\label{alfa1}
  (A,Y)\chi=(A,X)\gfi.
\end{equation}
By defining $X\alpha':=Y$ we have a collineation $\alpha'$ of $\PP_2$;
$\alpha:=(\rm{id}_{\call P_1},\alpha')$ is an automorphism of
$\PP_1\times \PP_2$.

Each line $\ell_2\in\call G_2$ such that\\
({\em a\/})~for every $X\in\call P_1$,
it holds
$( X \times\ell_2)\gfi=( X \times\ell_2)\alpha\chi$, and\\
({\em b\/})~the rank of $(\ach)\rist{\call P_1\times\ell_2}$ is greater than one,\\
will be called a {\em special line\/}.

The next result follows directly from  prop.~\ref{esistefi} and ({\em i\/}):

\begin{proposition}\label{speciale}
  Every line containing two points of $\call B$ is special.$\Box$
\end{proposition}

\begin{proposition}\label{eccezionali}
  Let $\ell_2$ be a special line,
  $X_1\in\call P_1\setminus A $,
  ${\bf A}_1:={\bf A}(\gfi)\cap(AX_1\times\ell_2)$, and
  ${\bf A}_2:={\bf A}(\ach)\cap(AX_1\times\ell_2)$.
  Then ${\bf A}_1={\bf A}_2$.
\end{proposition}

{\em Proof\/}.
We will use prop.~\ref{slm}.
Choose $B\in AX_1$ such that the dimension of
$U=( B \times \ell_2)\gfi=( B \times \ell_2)\ach$
is minimal.
If $U=\emptyset$, prop.~\ref{slm} implies
${\bf A}_1={\bf A}_2= B \times\ell_2$.
Otherwise, ${\bf A}_1$ and ${\bf A}_2$ have the same cardinality
$|{\bf A}_1|\leq2$.

Suppose $|{\bf A}_1|=1$. Take $X\in AB\setminus\{A,B\}$.
If $( A \times\ell_2)\gfi=( X \times\ell_2)\gfi$,
and $P$ is the unique point of $\ell_2$ satisfying
$\{(A,P)\gfi\}=( B \times\ell_2)\gfi$, then
${\bf A}_1=\{(B,P)\}$.
By definition of $\alpha$, $(A,P)\gfi=(A,P)\ach$, so
${\bf A}_2=\{(B,P)\}$.
If $( A \times\ell_2)\gfi\cap( X \times\ell_2)\gfi=\{(A,P_1)\gfi\}$,
we obtain
${\bf A}_1={\bf A}_2=\{(B,P_1)\}$ in a similar way.

Now suppose $|{\bf A}_1|=2$. Then there is a unique $C\in AB\setminus\{A,B\}$
such that
$( C \times\ell_2)\gfi=( C \times\ell_2)\ach$
has dimension zero.
Let $Q$ and $Q'$ be defined by
\[
  \{(A,Q)\gfi\}=( C \times\ell_2)\gfi,\hspace{2em}
  \{(A,Q')\gfi\}=( B \times\ell_2)\gfi,
\]
then
${\bf A}_1={\bf A}_2=\{(B,Q),(C,Q')\}$.$\Box$

\begin{proposition}\label{tesirido}
  If $\ell_2$ is a special line,
  $X_1\in\call P_1$ and $X_2\in\ell_2$, then\\
  $\{(X_1,X_2)\}\gfi=\{(X_1,X_2)\}\ach$.
\end{proposition}

{\em Proof\/}.
In case $X_1=A$ the statement holds by definition of $\alpha$.
From now on we assume $X_1\neq A$.
Let $r:=( A \times\ell_2)\gfi$,
$s:=( X_1 \times\ell_2)\gfi$.
Then $\dim r=1$, $\dim s\leq1$.
Since $\ell_2$ is special, it holds
$r=( A \times\ell_2)\ach$,
$s=( X_1 \times\ell_2)\ach$.
There are several cases to be considered:

{\em Case 1: $\dim s=-1$.\/} This is a trivial case.

{\em Case 2: $\dim s=0$.}\\
Since $s$ contains exactly one point, and
$\{(X_1,X_2)\}\gfi$, $\{(X_1,X_2)\}\ach$ are subsets of $s$,
the statement follows from prop.~\ref{eccezionali}.

{\em Case 3: $\dim s=1$, $r\cap s=\emptyset$.}\\
Let $A'\in AX_1\setminus\{A,X_1\}$.
Any two of the lines $r$, $s$ and $t:=( A' \times\ell_2)\gfi$ are skew.
Then, $(X_1,X_2)\gfi$ is the intersection of $s$ with the only line through
$(A,X_2)\gfi=(A,X_2)\ach$ intersecting $s$ and $t$.
Since the same property characterizes also $(X_1,X_2)\ach$, we have
$(X_1,X_2)\gfi=(X_1,X_2)\ach$.

{\em Case 4: $\dim s=1$, $r\cap s$ is a point $Z$.}\\
In this case, by props.~\ref{slm} and \ref{eccezionali},
${\bf A}_1={\bf A}_2$ either is empty or contains exactly one point.

{\em Case 4.1: $|{\bf A}_1|=1$.}\\
Let ${\bf A}_1=\{(A^*,P)\}$ ($A^*\in AX_1$),
then $A\neq A^*\neq X_1$, and
$( A^* \times\ell_2)\gfi=( A^* \times\ell_2)\ach$ is a point $Z'$.

The relations $r\vee s=r\vee\{Z'\}=s\vee\{Z'\}$ imply $Z'\not\in r\cup s$.
Then $(X_1,X_2)\gfi$ is the intersection point of $s$ and the line
containing both $Z'$ and $(A,X_2)\gfi$.
The same characterization holds for $(X_1,X_2)\ach$.

{\em Case 4.2: ${\bf A}_1=\emptyset$.}\\
Let $P,Q\in\ell_2$ defined by
$Z=(A,P)\gfi=(X_1,Q)\gfi$; then $P\neq Q$.

Let $g_X:=( X \times\ell_2)\gfi$ for $X\in AX_1$.
We now investigate on the following family of lines lying on the
plane $r\vee s$:
\[
  \call F:=\{g_X|X\in AX_1\}.
\]
First, $\call F$ is an injective family: indeed, if $X\neq X'$, then
$g_X\vee g_{X'}=r\vee s$, hence $g_X\neq g_{X'}$.

Now we prove that
{\em for every $T\in\ell_2$, the line
$h_T:=(AX_1\times T )\gfi$ belongs to $\call F$\/}.
Indeed, if $T=Q$, then $h_T=g_A=r$.
If $T\neq Q$, then $h_T$ contains $(X_1,T)\gfi$, that is a point not on
$r$.
Let $ Z'' :=r\cap h_T=(A,T)\gfi$.
There exists $C\in AX_1$ such that
$Z''=(C,Q)\gfi$.
It follows $g_C=(C,Q)\gfi\vee(C,T)\gfi=h_T$.

Next we prove that for every $X\in AX_1$,
{\em the line $g_X$ contains exactly one distinguished point $W_X$ which belongs
to no line of $\call F\setminus \{g_X\}$; if $U\in g_X\setminus W_X $,
then $U$ lies on exactly one line of
$\call F\setminus \{g_X\}$\/}.
\footnote{It should be noted that $\call F$ is a dual conic.}
Let $X'\in AX_1\setminus X $.
The lines $g_X$ and $g_{X'}$ intersect in exactly one point, say
$(X,Y)\gfi=(X',Y')\gfi$.
Then
\begin{equation}\label{efam}
  (AX_1\times Y' )\gfi=g_X,
\end{equation}
and every line $g_{X''}$, $X''\neq X$, meets $g_X$ in exactly one
point different from $W_X:=(X,Y')\gfi$; furthermore, distinct lines meet $g_X$ in
distinct points.
Conversely, by (\ref{efam}) every point of
$g_X\setminus W_X $ lies on a line $g_{X''}$ with
$X''\in AX_1\setminus X $.

The structure of $\call F$ allows the following considerations.
The distinguished points of the lines $r$ and $s$ are $W_A=(A,Q)\gfi$
and $W_{X_1}=(X_1,P)\gfi$, respectively.
The mapping \applicaz frs, defined by
$(A,X_2)\gfi f:=(X_1,X_2)\gfi$, can be characterized as follows:
$W_Af=Z$; $Zf=W_{X_1}$; for every $U\in r\setminus\{W_A,Z\}$,
$Uf$ is the intersection of $s$ with the unique line of $\call F\setminus\{r\}$
containing $U$.
Since for every $X\in AX_1$ we have
$g_X=( X \times\ell_2)\ach$, the same
geometric characterization holds for the mapping
$(A,X_2)\ach\mapsto(X_1,X_2)\ach$.

{\em Case 5: $r=s$.}\\
By property ({\em b\/}) for $\ell_2$, there is
$B\in\call P_1$ such that $( B \times\ell_2)\gfi\setminus r\neq\emptyset$.
Then, by prop.~\ref{slm}, there exists $B'\in AB$ such that
$( B' \times\ell_2)\gfi$ is a line $r'\neq r$.
If $X'_2\in\ell_2$, then
$(B',X'_2)\gfi=(B',X'_2)\ach$, because of the results of case 3 and 4.
Since $r'\neq s$, we can repeat the arguments in such cases replacing $A$
by $B'$ and $r$ by $r'$: indeed, the only properties of $A$ which are used
there are
$\gfi\rist{ A \times\ell_2}=\ach\rist{ A \times\ell_2}$ and
$\dim\left(( A \times\ell_2)\gfi\right)=1$.$\Box$

Now we can establish the main result of this section.

\begin{theorem}\label{main1}
  If conditions (i) and (ii), stated at the beginning of this section,
  hold, then there are a collineation $\alpha'$ of $\PP_2$ and a \ml\
  \applicaz{\varphi}{\overline{\PP}}{\PP'}\ such that, for
  $\alpha:=(\rm{id}_{\call P_1},\alpha')$, it holds
  $\gfi=\ach$.
\end{theorem}

{\em Proof\/}.
Take into account the \ml\ $\varphi$ of prop.~\ref{esistefi} and the
collineation $\alpha'$ defined by (\ref{alfa1}).
We will prove by induction on $t=0,1,\ldots,m$, that
\begin{equation}\label{hpindutt}
  \gfi\rist{\call P_1\times[A_0,A_1,\ldots,A_t]}=
  \ach\rist{\call P_1\times[A_0,A_1,\ldots,A_t]}.
\end{equation}
This is clear for $t=0$.
Next, assume that (\ref{hpindutt}) is true for some $t<m$.
Let $X_1\in\call P_1$, and
$X_2\in[A_0,A_1,\ldots,A_{t+1}]\setminus[A_0,A_1,\ldots,A_{t}]$.
Then there is $X_2'\in[A_0,A_1,\ldots,A_{t}]$ such that
$X_2\in A_{t+1}X_2'=:\ell_2$.

If $X\in\call P_1$, then by the induction hypothesis
$\{(X,X_2')\}\gfi=\{(X,X_2')\}\ach$, and, by prop.~\ref{esistefi},
$\{(X,A_{t+1})\}\gfi=\{(X,A_{t+1})\}\ach$.
It follows $( X \times\ell_2)\gfi=( X \times\ell_2)\ach$.
This implies that $\ell_2$ is a special line (see also assumption ({\em ii\/})).
Our statement follows then from prop.~\ref{tesirido}.$\Box$


\section{The case $\rad 2\chi\neq\emptyset$}

For each $P\in\call P_1$, let \applicaz{\chi_P}{\PP_2}{\PP'}\ denote the \ml\
which is defined by $\{X\} \chi_P=\{(P,X)\}\chi$.

\begin{theorem}\label{degenere}
  Assume that there exists $A\in\call P_1$ such that
  $\dim{\bf A}(\chi_A)\leq m-2$, and,
  for every $P\in\call P_1$,
  ${\bf A}(\chi_A)\subset{\bf A}(\chi_P)$.
  Also assume that there is a
  basis
  $\{A_0,A_1,\ldots,A_p\}$ of a subspace $D_2$,
  which is complementary with ${\bf A}(\chi_A)$ in $\PP_2$,  and a plane $\call E\subset\call P_1$
  such that for every $i=1,\ldots,p$,
  $\chi\rist{\call E\times\{A_0,A_i\}}$ is global and of rank greater than two.
  Then there exist a \ml\ \applicaz{\varphi}{\overline{\PP}}{\PP'}\ and
  a collineation $\alpha'$ of $\PP_2$ such that, for
  $\alpha:=(\rm{id}_{\call P_1},\alpha')$, it holds
  $\gfi=\ach$.
\end{theorem}

{\em Proof\/}.
There is a subspace $D_1\subset\call P_1$ such that $D_1$ and $\rad 1\chi$ are
complementary subspaces of $\PP_1$ and $\call E\subset D_1$.
Next, $\rad 2\chi={\bf A}(\chi_A)$, and $D_2\subset{\bf D}(\chi_A)$.
Then theor.~\ref{main1} applies to the \ml\ $\chi'$ that is defined in
prop.~\ref{nondegen}.
So, there are a collineation $\beta'$
of $D_2$ and a \ml\
\applicaz{\varphi'}{[(D_1\times D_2)\gamma]}{\PP'}\ such that, for
$\beta:=({\rm id}_{D_1},\beta')$, it holds
$\beta\chi'=(\gamma\varphi')\rist{D_1\times D_2}$.

It is possible to extend $\beta'$ to a collineation $\alpha'$ of $\PP_2$
such that $(\rad 2\chi)\alpha'=\rad 2\chi$.
Since $\gamma$ is regular,
\[
  U:=\left[((\rad 1\chi)\times\call P_2)\gamma\right]\vee
     \left[(\call P_1\times(\rad 2\chi))\gamma\right]
\]
is a subspace of $\overline{\PP}$ which is complementary with
$[(D_1\times D_2)\gamma]$.
Let \applicaz{\varphi}{\overline{\PP}}{\PP'}\
be a \ml\ extending
$\varphi'$ such that $U\varphi=\emptyset$.
This $\varphi$ satisfies the required property.$\Box$

\vspace{6mm}

\noindent
Alessandro Bichara\\
Dipartimento di Metodi e Modelli Matematici\\
Universit\`a ``La Sapienza''\\
via Scarpa 16, I-00161 Roma, Italy\\
bichara@dmmm.uniroma1.it

\vspace{3mm}
\noindent
Hans Havlicek\\
Technische Universit\"at Wien\\
Abteilung f\"ur Lineare Algebra und Geometrie\\
Wiedner Hauptstrasse 8-10/1133, A-1040 Wien, Austria\\
havlicek@geometrie.tuwien.ac.at

\vspace{3mm}
\noindent
Corrado Zanella\\
Dipartimento di Matematica Pura ed Applicata\\
Universit\`a di Padova\\
via Belzoni 7, I-35131 Padova, Italy\\
zanella@math.unipd.it

\end{document}